\documentclass[graybox]{svmult}

\usepackage{type1cm}        % activate if the above 3 fonts are
\usepackage{makeidx}         % allows index generation
\usepackage{graphicx}        % standard LaTeX graphics tool
\usepackage{multicol}        % used for the two-column index
\usepackage[bottom]{footmisc}% places footnotes at page bottom

\usepackage{newtxtext}       % 
\usepackage{newtxmath}       % selects Times Roman as basic font
\usepackage{amssymb} % ajoute par thomas, a retirer ensuite
\usepackage[utf8]{inputenc}

\usepackage[all,cmtip]{xy}
\usepackage{amsmath}
\usepackage{xcolor}

\def\c1{\operatorname{c_1}}
\def\c2{\operatorname{c_2}}

\def\ZZ{{\mathbb Z}}

\def\PP{{\mathbb P}}

\def\O{{\mathcal O}}
\def\I{{\mathcal J}}

\def\c{\mathfrak{c}}

\def\cong{\simeq}

\def\+{\oplus}               % direct sum
\def\*{\otimes}                  % tensor product
\def\Pic{\operatorname{Pic}}

\makeindex             % used for the subject index
\begin{document}

\title*{A note on Severi varieties of nodal curves on Enriques surfaces}
 \titlerunning{Severi varieties on Enriques surfaces} 
\author{Ciro Ciliberto, Thomas Dedieu, Concettina Galati and Andreas Leopold Knutsen}
 \authorrunning{C.~Ciliberto, T.~Dedieu, C.~Galati and A.~L.~Knutsen}
\institute{Ciro Ciliberto \at Dipartimento di Matematica, Universit\`a
  di Roma Tor Vergata, Via della Ricerca Scientifica, 00133 Roma,
  Italy, \email{cilibert@axp.mat.uniroma2.it} 
\and Thomas Dedieu \at 
Institut de Mathématiques de Toulouse, UMR5219.
Université de Toulouse, CNRS.
UPS IMT, F-31062 Toulouse Cedex 9, France.
\email{thomas.dedieu@math.univ-toulouse.fr}
\and
Concettina Galati \at Dipartimento di Matematica e Informatica,
Universit\`a della Calabria, via P. Bucci, cubo 31B, 87036 Arcavacata di
Rende (CS), Italy, \email{galati@mat.unical.it} 
\and Andreas Leopold Knutsen \at Department of Mathematics, University
of Bergen, Postboks 7800, 
5020 Bergen, Norway \email{andreas.knutsen@math.uib.no}}
\maketitle

\abstract*{Each chapter should be preceded by an abstract (no more than 200 words) that summarizes the content. The abstract will appear \textit{online} at \url{www.SpringerLink.com} and be available with unrestricted access. This allows unregistered users to read the abstract as a teaser for the complete chapter.
Please use the 'starred' version of the \texttt{abstract} command for typesetting the text of the online abstracts (cf. source file of this chapter template \texttt{abstract}) and include them with the source files of your manuscript. Use the plain \texttt{abstract} command if the abstract is also to appear in the printed version of the book.}

\abstract{Let $|L|$ be a linear system on a smooth complex Enriques
surface $S$ whose general member is a smooth and irreducible curve of genus
$p$, with $L^ 2>0$, and let $V_{|L|, \delta} (S)$ be the Severi
variety of irreducible $\delta$-nodal curves in $|L|$. We denote by
$\pi:X\to S$ the universal covering of $S$. In this note we compute
the dimensions of the irreducible components $V$ of $V_{|L|, \delta}
(S)$.
In particular we prove that, if $C$ is the curve corresponding to a
general element $[C]$ of $V$, then the codimension of $V$ in $|L|$ is
$\delta$ if $\pi^{-1}(C)$ is irreducible in $X$ and it is $\delta-1$
if $\pi^ {-1}(C)$ consists of two irreducible components.}

\section{Introduction}
Let $S$  be a smooth complex projective  surface and $L$ a line bundle on $S$ such that the complete linear system $|L|$ contains smooth, irreducible curves (such a line bundle, or linear system, is often called a {\it Bertini system}). Let
\[
p:=p_a(L)=  \frac{1}{2} L \cdot (L+K_S) +1,
\]
be the arithmetic genus of any curve in  $|L|$. 

For any integer $0 \leq \delta \leq p$, consider the locally closed, functorially defined subscheme of $|L|$  
\[
V_{|L|, \delta} (S) \; \; \mbox{or simply} \; \; V_{|L|, \delta} 
\] 
parameterizing irreducible curves in $|L|$ 
having only $\delta$ nodes as singularities; this is  called the {\it Severi variety} of $\delta$-nodal curves in $|L|$. We will let $g:=p-\delta$, the geometric genus of the curves in $V_{|L|, \delta}$.

It is well-known that, if $ V_{|L|, \delta}$ is non-empty, then
all of its irreducible components $V$ have dimension $\dim(V) \geq
\dim |L|-\delta$. More precisely, the Zariski tangent space to
$V_{|L|,\delta}$ at the point corresponding to $C$ is
\begin{equation}\label{eq:Sevvar}
T_{[C]} V_{|L|, \delta}\simeq H^0(L \* \I_{N})/<C>,
\end{equation}
where  $ \I_{N}= \I_{N|S}$ is the ideal sheaf of subscheme $N$ of $S$ consisting of the $\delta$ nodes of $C$ (see, e.g., \cite[\S 1]{CS}). Thus, $V_{|L|, \delta}$ is {\it smooth of dimension $\dim|L|-\delta$} at $[C]$ if and only if the set of nodes $N$ imposes independent conditions on $|L|$. In this case, $V_{|L|, \delta}$ is said to be {\it regular} at $[C]$.  An irreducible component $V$ of $V_{|L|, \delta}$ will be said to be  {\it regular} if the condition of regularity is satisfied at any of its points, equivalently, if it is smooth of dimension $\dim |L|-\delta$.

The {\it existence and regularity problems of $V_{|L|, \delta} (S)$}
have been studied in many cases and are the most basic problems one
may ask on Severi varieties.  We only mention some of known
results. In the case $S\simeq\mathbb P^2$ , Severi 
proved the existence and regularity of $V_{|L|, \delta} (S)$ 
in \cite{Sev}. The description of the tangent
space is due to Severi and later to Zariski \cite{Zar}. The existence
and regularity of $V_{|L|, \delta} (S)$ when $S$ is of general type
has been studied in \cite{CS} and \cite{CC}. Further regularity
results are provided in \cite{F1}. More recently Severi varieties on
K3 surfaces have received a lot of attention for many reasons. In this
case Severi varieties are known to be regular (cf. \cite{Tan}) and are
nonempty on general K3 surfaces by Mumford and Chen
(cf. \cite{MM}, \cite{Ch}).

As far as we know, Severi varieties on Enriques surfaces have not been
studied yet, apart from \cite[Thm.~4.12]{DS} which limits the
singularities of a general member of the Severi variety $V_{|L|}^g$ of
irreducible genus $g$ curves in $|L|$, and gives a sufficient
condition for the density of the latter in the Severi variety
$V_{|L|,p-g}$ of $(p-g)$-nodal curves.
In particular, the existence problem is mainly open and
we intend to treat it in a forthcoming article. The result of this
paper is Proposition~\ref{prop:reg-enr}, which answers
the regularity question for Severi varieties of nodal curves on
Enriques surfaces.

\section{Regularity of Severi varieties on Enriques surfaces}
\label{S:severisetup}
Let $S$ be a smooth Enriques surface, i.e. a smooth
complex surface with nontrivial canonical bundle $\omega_S \not \cong
\mathcal O_S$, such that $\omega_S^{\* 2} \cong \O_S$ and $H^1(\mathcal O_S)=0$. 
We denote linear (resp.\ numerical) equivalence by $\sim$ (resp.\
$\equiv$). 

Let
$L$ be a line bundle on $S$ such that $L^ 2>0$. It is well-known that
$|L|$ contains smooth, irreducible curves if and only if it contains
irreducible curves (see \cite[Thm. 4.1 and Prop. 8.2]{cos1}); in other
words, {\it on Enriques surfaces the Bertini linear systems are the
linear systems that contain irreducible curves}. Moreover, by
\cite[Prop. 2.4]{cos2}, this is equivalent to $L$ being nef and not of
the form $L \sim P+R$, with $|P|$ an elliptic pencil and $R$ a smooth
rational curve such that $P \cdot R=2$ (in which case $p=2$). If $|L|$
is a Bertini linear system, the adjunction formula, the Riemann--Roch
theorem,
and Mumford vanishing yield that
\[ L^2 =2(p-1) \; \; \mbox{and} \; \; \dim |L|=p-1\]
(see, e.g., \cite{cos1,cd}).

Let $K_S$ be the canonical divisor. It defines an \'etale double cover
\begin{equation} \label{eq:dc}
 \pi: X \longrightarrow S 
\end{equation}
where $X$ is a smooth, projective $K3$ surface (that is, $\omega_X \cong \mathcal O_X$ and $H^1(\mathcal O_X)=0$), endowed with a fixed-point-free involution $\iota$, which is the universal covering of $S$. Conversely, the quotient of any $K3$ surface by a fixed-point-free involution is an Enriques surface. 

\begin{figure}[htbp] 
\noindent
%\hfill\hfill\hfill
\begin{minipage}[t]{0.47\textwidth}%\label{figura2}
\begin{center}
\includegraphics[height=5 cm, width=4cm]{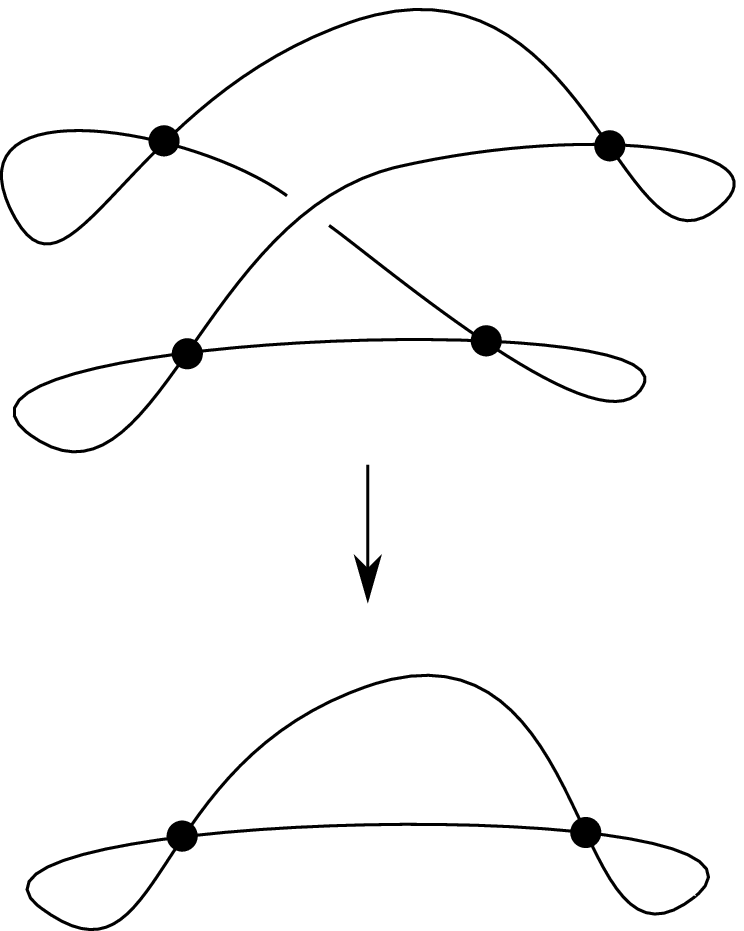}
\end{center}
%{fig2.pdf} 
\sidecaption
\caption{$\eta _{\tilde C} = \nu_C^*(\eta_C) \neq 0$}
%{Irreducible $\pi^{-1}C$}
\label{fig:irr}
\end{minipage}\hfill
\begin{minipage}[t]{0.50\textwidth}%\label{figura1}
\begin{center}
\includegraphics[height=5 cm, width=4cm ]{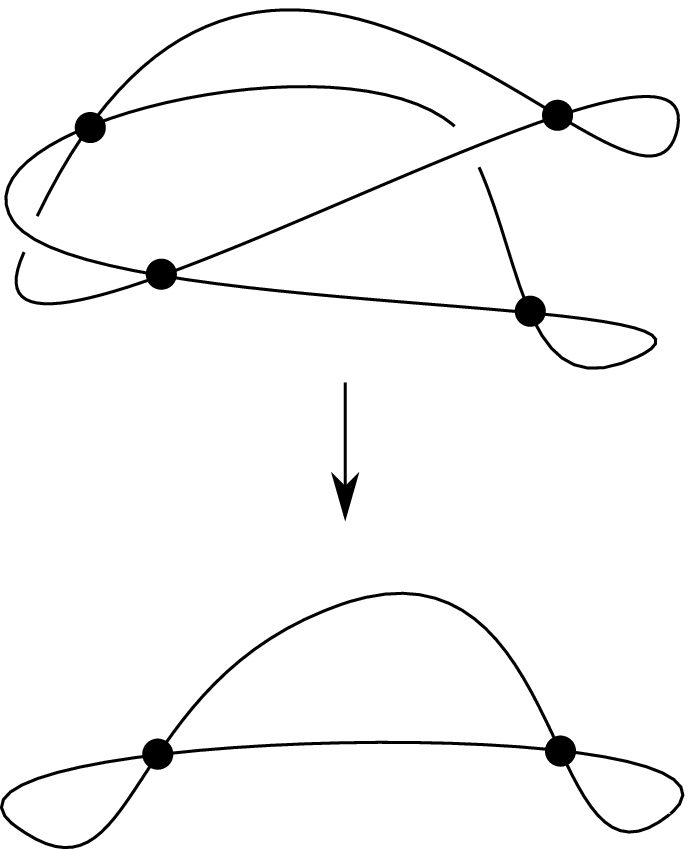}
%{fig1.pdf}
\end{center}
\sidecaption
\caption{$\eta _{\tilde C} =\nu_C^*(\eta_C) = 0$}
%{Reducible $\pi^{-1}C$}
\label{fig:red}
\end{minipage}
\end{figure}

Let $C \subset S$ be a reduced and irreducible curve of genus $g \geq
2$. We will henceforth denote by $\nu_C:\widetilde{C} \to C$ the
normalization of $C$ and define $\eta_C:=\mathcal O_C(K_S)=\mathcal O_C(-K_S)$, a
nontrivial $2$-torsion element in $\Pic^0 C$, and
$\eta_{\widetilde{C}}:=\nu_C^*\eta_C$. The fact that $\eta_C$ is
nontrivial follows from the cohomology of the restriction sequence
\[ \xymatrix{
0 \ar[r] & \mathcal O_S(K_S-C)  \ar[r] & \mathcal O_S(K_S) \ar[r] & \eta_C \ar[r] & 0,}
\] which yields $h^0(\eta_C)=h^1(K_S-C)=h^1(C)=0$, the latter vanishing as $C$ is big and nef.
One has the fiber product 
\[ \xymatrix{
(\pi^{-1}C) \times_C \widetilde{C} \ar[r] \ar[d] & \widetilde{C} \ar[d]^{\nu_C} \\
(\pi^{-1}C) \ar[r]_{\pi_{|_{\pi^{-1}(C)}}} & C,}
\]
where $\pi_{|_{\pi^{-1}(C)}}$ and the upper horizontal map are the double coverings induced respectively by $\eta_{C}$ and $\eta_{\widetilde{C}}$. 
By standard results on coverings of complex manifolds (cf. \cite[Sect. I.17]{BPHV}), two cases may happen:
\begin{itemize}
\item $\eta_{\widetilde{C}} \not \cong \mathcal O_{\widetilde{C}}$ and
$\pi^{-1}C$ is irreducible, as in Fig.~\ref{fig:irr};
\item $\eta_{\widetilde{C}} \cong \mathcal O_{\widetilde{C}}$ and $\pi^{-1}C$
consists of two irreducible components conjugated by the involution
$\iota$. These two components are {\it not} isomorphic to $C$, as
$\eta_C$ is nontrivial, as in Fig.~\ref{fig:red}
(each component of $\tilde C$ is a partial normalization of $C$).
\end{itemize}

As mentioned in the Introduction, it is well-known that {\it any}
irreducible 
component of a Severi variety on a $K3$ surface is regular when
nonempty (see, e.g., \cite[Ex.~1.3]{CS}; 
see also \cite[\S 4.2]{DS}).
The corresponding result on Enriques surfaces is the following.

First note that, in the above notation, the dimension of the Severi
variety of genus $g=p_g(C)$ curves in $|L|=|C|$ at the point $[C]$
satisfies the inequality
\begin{equation}
\label{ineq:equigen-tg}
\dim _{[C]} \bigl( V_{|L|} ^g \bigr)
\geq h^0(\omega _{\tilde C} \otimes \eta _{\tilde C})
= \begin{cases}
g-1 & \text{if } \eta_{\tilde C} \not\simeq \mathcal O_{\tilde C} \\
g & \text{if } \eta_{\tilde C} \simeq \mathcal O_{\tilde C}
\end{cases}
\end{equation}
(see \cite[Proofs of Thm.~4.12 and Cor.~2.7]{DS}).
Our result implies that the latter is in fact an equality when $C$ is
nodal, and gives a concrete geometric description of the situation in
both cases.

\begin{proposition} \label{prop:reg-enr}
 Let $L$ be a Bertini linear system, with $L^2>0$, on a smooth
Enriques surface $S$.  Then the Severi variety $V_{|L|, \delta}(S)$ is
smooth and every irreducible component $V\subseteq V_{|L|, \delta}(S)$
has either dimension $g-1$ or $g$; in the former case the component is
regular. Furthermore, with the notation introduced above,
\begin{enumerate}
\item for any curve $C$ in a $(g-1)$-dimensional irreducible
component $V$, $\pi^{-1}C$ is irreducible
(whence an element in $V_{|\pi^*L|,2\delta}(X)$); 
\item for any $g$-dimensional component $V$, there is a line
bundle $L'$ on $X$ with $(L')^2=2(p-d)-2$
 and $L' \cdot \iota^*L'=2d$
for some integer $d$ satisfying 
\[ \frac{p-1}{2} \leq d \leq \delta,\]
such that 
$\pi^*L \cong L' \* \iota^*L'$,
and the curves parametrized by 
$V \subseteq  V_{|L|, \delta}(S)$ 
are the birational images by $\pi$
of the curves in $V_{|L'|, \delta-d}(X)$ intersecting their conjugates
by $\iota$ transversely (in $2d$ points). In other words, for any $[C]
\in V$, we have $\pi^{-1}C=Y+\iota(Y)$, with $[Y] \in V_{|L'|,
\delta-d}(X)$ and $[\iota(Y)] \in V_{|\iota^*L'|, \delta-d}(X)$
intersecting transversely.

Furthermore, if $L' \cong \iota^*L'$, which is the case if $S$ is
general in moduli, then $d= \frac{p-1}{2}$ and $L \sim 2M$, for some
$M \in \Pic S$ such that $M^2=d$.
\end{enumerate}
\end{proposition}

We will henceforth refer to components of dimension $g-1$ as {\it
regular} and the ones of dimension $g$ as {\it nonregular}.
Note however that from a parametric perspective the Severi variety
has the expected dimension and is smooth in both cases, as the fact that
\eqref{ineq:equigen-tg} is an equality indicates; we do not dwell on
this here, and refer to \cite{DS} for a discussion of the differences
between the parametric and Cartesian points of view
(the latter is the one we adopted in this text).

Note that Proposition~\ref{prop:reg-enr} does not assert that the
Severi variety $V_{|L|,\delta}$ is necessarily non-empty: in such a
situation, $V_{|L|,\delta}$ does not have any irreducible component
and the statement is empty.

\begin{proof}
Pick any curve $C$ in an irreducible component $V$ of $V_{|L|, \delta}(S)$ . Let $f: \widetilde{S} \to S$ be the blow-up of $S$ at $N$, the scheme of the $\delta$ nodes of $C$, denote by $\mathfrak{e}$ the (total) exceptional divisor and by $\widetilde{C}$  the strict transform of $C$. Thus $f_{|\widetilde{C}}=\nu_C$ and we have
\[ K_{\widetilde{S}} \sim f^*K_S+\mathfrak{e} \; \; \mbox{and} \; \; \widetilde{C} \sim f^*C-2\mathfrak{e}.\]
From the restriction sequence
\[
\xymatrix{
0 \ar[r] & \mathcal O_{\widetilde{S}}(\mathfrak{e}) \ar[r] & \mathcal O_{\widetilde{S}}(\widetilde{C}+\mathfrak{e}) \ar[r] & \omega_{\widetilde{C}}(\eta_{\widetilde{C}}) \ar[r] & 0 }\]
we find 
\begin{eqnarray}
  \dim T_{[C]} V_{|L|, \delta}(S) = \dim |L \* \I_N| & = & h^0(L \* \I_N)-1=h^0(f^*L-\mathfrak{e})-1 \nonumber \\ 
& = & h^0(\mathcal O_{\widetilde{S}}(\widetilde{C}+\mathfrak{e}))-1 = h^0(\omega_{\widetilde{C}}(\eta_{\widetilde{C}}))\nonumber \\
& = & \begin{cases}
g-1, & \mbox{if} \; \eta_{\widetilde{C}} \not \cong \mathcal O_{\widetilde{C}}, \\
g, & \mbox{if} \; \eta_{\widetilde{C}} \cong \mathcal O_{\widetilde{C}}. \label{nonregular}
\end{cases}
\end{eqnarray}

In the upper case, by \eqref{eq:Sevvar}, we have that $V_{|L|, \delta}$ is smooth at $[C]$ of dimension $g-1=p-\delta-1=\dim |L \* \I_N|$. 

Assume next that we are in the lower case. Then, by the discussion
prior to the proposition, we have $\pi^{-1}C=Y+\iota(Y)$ for
an irreducible curve $Y$ on $X$, such that $\pi$ maps both $Y$ and
$\iota(Y)$ birationally, but not isomorphically, to $C$. In
particular, $Y$ and $\iota(Y)$ have geometric genus
$p_g(Y)=p_g(\iota(Y))=p_g(C)=p-\delta=g$.  Set $L':=\mathcal O_X(Y)$ and
$2d:=Y \cdot \iota(Y)$. Note that $d$ is an integer because, if
$y=\iota(x)\in Y\cap \iota(Y)$, then $\iota(y)=x\in Y\cap
\iota(Y)$. Since $Y \cong \iota(Y)$ and $\pi$ is \'etale, both $Y$ and
$\iota(Y)$ are nodal with $\delta-d$ nodes and they intersect
transversely at $2d$ points, which are pairwise conjugate by $\iota$,
and therefore map to $d$ nodes of $C$. Hence $d \leq \delta$. We have
\begin{equation} \label{eq:genY}
p_a(Y)=p_a(\iota(Y))=g+\delta-d=p-\delta+\delta-d=p-d.
\end{equation}
whence
\begin{eqnarray*} 
(L')^2&=&2(p-1-d).
\end{eqnarray*}
By the Hodge index theorem, we have
\[ 4(p-1-d)^2=\left((L')^2\right)^2 = (L')^2 (\iota^*L')^2 \leq \left(L' \cdot \iota^*L'\right)^2=4d^2,\]
whence $p-1\leq 2d$. 

By the regularity of Severi varieties on $K3$ surfaces, any irreducible component of  $V_{|L'|, \delta-d}(X)$ has dimension $\dim |L'|-(\delta-d)=p_g(Y)=g$. Hence, $V$ is $g$-dimensional; more precisely, the curves parameterized  by $V$ are the (birational) images by $\pi$ of the curves in an irreducible component of $V_{|L'|, \delta-d}(X)$ intersecting their conjugates by $\iota$ transversely (in $2d$ points). By \eqref{nonregular}, it also follows that $\dim V=\dim T_{[C]} V_{|L|, \delta}(S)$, so that
$[C]$ is a smooth point of $V_{|L|, \delta}(S)$. 

To prove the final assertion of the proposition, observe that, by the
regularity of Severi varieties on $K3$ surfaces, we may deform $Y$ and
$\iota(Y)$ on $X$ to irreducible curves $Y'$ and $\iota(Y')$ with any
number of nodes $\leq \delta-d$ and intersecting transversally in $2d$
points; in particular, we may deform $Y$ and $\iota(Y)$ to {\it
smooth} curves $Y'$ and $\iota(Y')$.  Thus, $C':=\pi(Y')$
is a member of $V_{|L|,d}$, whence of geometric genus $p-d$. Since
$\dim |Y'|=p_a(Y')=p_g(C')=p_a(C')-d=p-d$, the component of
$V_{|L|,d}$ containing $[C']$ has dimension $\dim |L|-d+1=p-d$.
We thus have $\dim |L \* \I_{N'}|=\dim
|L|-d+1$, where $N'$ is the set of $d$ nodes of $C'$, hence $N'$ does
not impose independent conditions on $|L|$.

Assume now that $L' \cong \iota^*L'$, which ---~as is well-known (see,
e.g., \cite[\S 11]{DK})~--- is the case occurring for generic $S$, as
then $\Pic X$ is precisely the invariant part under $\iota$ of
$H_2(X,\ZZ)$. Then $2d= L' \cdot \iota^*L'=(L')^2=2(p-1-d)$, so that
$p-1=2d$.  Since
$L^2=2(p-1)=4d$ and $N'$ does not impose independent conditions on
$|L|$, by \cite[Prop. 3.7]{Kn} there is an effective divisor $D
\subset S$ containing $N'$ satisfying $L-2D \geq 0$ and
\begin{equation}
\label{eq:man} 
L\cdot D-d \leq D^2 \stackrel{(i)}{\leq} \frac{1}{2}L \cdot D
\stackrel{(ii)}{\leq} d,
\end{equation}
with equality in (i) or (ii) only if $L \equiv 2D$;
moreover, since  $L-2D \geq 0$, the numerical equivalence
$L \equiv 2D$ implies the linear equivalence $L \sim 2D$.
Now since $N' \subset D$, we must have $L
\cdot D =C' \cdot D \geq 2d$, hence the inequalities in \eqref{eq:man}
are all equalities, and thus $D^2=d$ and $L \sim 2D$.
\end{proof}

The following corollary is a straightforward consequence of Prop. \ref{prop:reg-enr} and the fact that the nodes on curves in a regular component in a Severi variety (on any surface and in particular on a K3 surface) can be independently smoothened.
\begin{corollary} \label{cor:reg-enr}
 If a Severi variety $V_{|L|, \delta}$ on an Enriques surface has a regular (resp., nonregular) component, then for any $0 \leq \delta' \leq \delta$ (resp.,
$d \leq \delta' \leq \delta$, with $d$ as in Prop. \ref{prop:reg-enr}), also $V_{|L|, \delta'}$ contains a 
regular (resp., nonregular) component.
\end{corollary}

\begin{acknowledgement}
The first and third author want to thank the organisers of BGMS Indam Workshop for the opportunity to give a talk. 
The first author acknowledges the MIUR Excellence Department Project
awarded to the Department of Mathematics, University of Rome Tor
Vergata, CUP E83C18000100006. 
The first and second authors were members of project FOSICAV, which 
has received funding from the European Union's Horizon
2020 research and innovation program under the Marie
Sk{\l}odowska-Curie grant agreement n.~652782.
The third author has been partially supported by GNSAGA of INDAM.
The last author has been partially supported by grant n.~261756 of the Research
Council of Norway and by Bergen Research Foundation.

\end{acknowledgement}

\end{document}